\documentclass[12pt,english,reqno]{amsart}
\usepackage{amssymb,amsmath,amstext,amsthm,amsfonts}
\usepackage[dvips]{graphicx}

\newtheorem{theorem}{Theorem}[section]
\newtheorem{lemma}[theorem]{Lemma}

\newtheorem{corollary}[theorem]{Corollary}

\def \RR {{\mathbb R}}
\def \ZZ {{\mathbb Z}}
\def \NN {{\mathbb N}}

    \def \co {\mathcal{O}}

\newcommand{\dem}{\begin{proof}}
\newcommand{\cqd}{\end{proof}}

\newcommand{\leb}{\operatorname{Leb}}

\newcommand{\qand}{\quad\mbox{and}\quad}

\begin{document}

\title[Backward volume contraction for endomorphisms]{Backward volume contraction for endomorphisms  with eventual volume
expansion}

\author{Jos\'e F. Alves}
\address{Departamento de Matem\'atica Pura, Faculdade de Ci\^encias do Porto\\
Rua do Campo Alegre 687, 4169-007 Porto, Portugal}
\email{jfalves@fc.up.pt}

\author{Vilton  Pinheiro}
\address{Departamento de Matem\'atica, Universidade Federal da Bahia\\
Av. Ademar de Barros s/n, 40170-110 Salvador, Brazil.}
\email{viltonj@ufba.br}

\author{Armando Castro}
\address{Departamento de Matem\'atica, Universidade Federal da Bahia\\
Av. Ademar de Barros s/n, 40170-110 Salvador, Brazil.}
\email{armando@im.ufba.br}

\date{\today}

\thanks{Work carried out at the  Federal University of
Bahia. Partially supported by FCT through CMUP and UFBA}

\maketitle

\begin{abstract}
We consider smooth maps on compact Riemannian manifolds. We prove
that under some mild condition of eventual volume expansion
Lebesgue almost everywhere  we have uniform backward volume
contraction  on every pre-orbit for Lebesgue almost every point.

\end{abstract}

\section{Statement of results}

Let $M$ be a compact Riemannian manifold and let $\leb$ be a
volume form on $M$ that we call Lebesgue measure. We take $f\colon
M\to M$ any smooth map.
 Let $0<a_1\le a_2 \le a_3\le \dots$ be a sequence  converging
to infinity. We define
 \begin{equation}\label{aga}
    h(x)=\min\{n>0 \colon |\det Df^n(x)|\ge a_n\},
\end{equation}
if this minimum exists, and $h(x)=\infty$, otherwise. For $n\ge
1$, we take
 \begin{equation}\label{gamma}
   \Gamma_n=\{x\in M \colon h(x) \ge n\}.
\end{equation}

\begin{theorem} \label{JCtheo} Assume that $h\in L^p(\leb)$, for some $p>3$, and
take $\gamma<(p-3)/(p-1)$. Choose any sequence $0<b_1\le b_2 \le
b_3\le \dots$ such that $b_kb_n\ge b_{k+n}$ for every $k,n\in
\NN$, and assume that there is $n_0\in\NN$ such that $b_n\le
\min\left\{a_n,\leb(\Gamma_n)^{-\gamma}\right\} $ for every $n\ge
n_0$. Then, for $\leb$ almost every $x\in M$, there exists $C_x>0$
such that $| \det Df^n(y)|>C_x b_n$ for every $y\in f^{-n}(x).$
 \end{theorem}

\smallskip

We say that $f\colon M\to M $ is {\em eventually volume expanding}
if there exists $\lambda>0$ such that for Lebesgue almost every
$x\in M$
 \begin{equation}
 \sup_{n\ge 1}\frac1n\log|\det Df^n(x)|> \lambda.
 \end{equation}
Let $h$ and $\Gamma_n$ be defined as in~(\ref{aga})
and~(\ref{gamma}), associated to the sequence $a_n=e^{\lambda n}$.

\medskip

 \begin{corollary} \label{Mtheo} If $f$ is eventually volume expanding, then for Le\-bes\-gue almost
every point $x\in M$ there are $C_x>0$ and $\sigma_n\to \infty$
such that $| \det Df^n(y)|>C_x \sigma_n $ for every $y\in
f^{-n}(x)$. Moreover, given $\alpha>0$ there is $\beta>0$ such
that
\begin{enumerate}
\item   if $\leb(\Gamma_n)\le \co(e^{-\alpha n})$, then we may
take $ \sigma_n \ge e^{\beta n}$;
    \item if
$\leb(\Gamma_n)\le \co(e^{-\alpha n^\tau})$ for some $\tau>0$,
then we may take $\sigma_n \ge e^{\beta n^\tau}$;
    \item if
$\leb(\Gamma_n)\le \co(n^{-\alpha})$ and $\alpha>2$, then we may
take $\sigma_n \ge n^\beta$.
\end{enumerate}
 \end{corollary}

 \medskip

Specific rates  will be obtained in Section~\ref{se.examples} for
some eventually volume expanding endomorphisms. In particular,
non-uniformly expanding maps such as quadratic maps and Viana maps
will be considered.

%

\section{Concatenated collections}

 Let  $(U_n)_n$ be a collection of measurable subsets
 of $M$ whose union covers a full Lebesgue measure subset of  $M$. We say that
$(U_n)_n$ is a {\em concatenated collection}
  if:
  $$x\in U_n \qand f^n(x)\in U_m\quad\Rightarrow\quad x\in
  U_{n+m}.$$
Given $x\in \bigcup_{n\ge 1} U_n$, we define
 $u(x)$ as the minimum
$n\in\NN$ for which $x\in U_n$. Note that by definition we have
$x\in U_{u(x)}$. 
We define  the {\em chain generated by $x\in \bigcup_{n\ge 1}
U_n$} as
 $C(x)=\{x,f(x),\dots,f^{u(x)-1}(x)\}.$

\begin{lemma}
\label{JClemma1} Let $(U_n)_n$ be a
 concatenated collection. If
$$\sum_{n\ge 1}\sum_{j=0}^{n-1}\leb(f^j(u^{-1}(n)))<\infty,$$ then
we have $ \sup\left\{\,u(y)\ \colon \; y\in \bigcup_{n\ge 1}
U_n\;\mbox{and}\,\; x\in C(y)\,\right\}<\infty $ for Leb\-es\-gue
almost every $x\in M$.
\end{lemma}

\dem
 Assume that for a given $x\in M$ there exists an infinite
number of chains
 $C_j=\left\{y_j,f(y_j), \dots ,f^{s_j-1}(y_j)\right\}$, $ j\ge 1$, containing $x$
 with $s_j\to\infty$. For each $j\ge1$ let
$1\le r_j<s_j$ be such that $x=f^{r_j}(y_j)$.
 First we verify that
$\lim r_j=\infty$. If not, then replacing by a subsequence, we may
assume that there is $N>0$ such that $r_j<N$ for every $j\ge1$.
This implies that $y_j\in\bigcup_{i=1}^{N}f^{-i}(x)$ for every
$j\ge1$. Since $\#(\bigcup_{i=1}^{N}f^{-i}(x))<\infty$ and the
number of chains is infinite, we have a contradiction.
 Since
$r_j\to\infty$ and $x=f^{r_j}(y_j)\in f^{r_j}(u^{-1}(s_j))$, then
we have $x\in\bigcup_{n\ge k}\bigcup_{j=0}^{n-1}f^j(u^{-1}(n))$
for every $k\ge 1$. Since we are assuming $\sum_{n\ge
1}\sum_{j=0}^{n-1}\leb(f^j(u^{-1}(n)))<\infty$, we have
$\leb\big(\bigcup_{n\ge
k}\bigcup_{j=0}^{n-1}f^j(u^{-1}(n))\big)\to 0, $ when
$k\to\infty$. This completes the proof of Lemma~\ref{JClemma1}.
\cqd

\begin{lemma}
\label{JClemma2} Let $(U_n)_n$ be a concatenated collection. If
$$\sup\left\{\,u(y)\ \colon \; y\in \cup_{n\ge 1}
U_n\;\mbox{and}\,\; x\in C(y)\,\right\}\le N,$$ 
then $f^{-n}(x)\subset U_{n}\cup \dots \cup U_{n+N}$ for all
$n\ge1$.
\end{lemma}

\dem Assume that $\sup\left\{\,u(y)\ \colon \; y\in \cup_{n\ge 1}
U_n\;\mbox{and}\,\; x\in C(y)\,\right\}\le N$, and take $z\in
f^{-n}(x)$.  Let $z_j=f^j(z)$ for each $j\ge 0$. We distinguish
the cases $x\in C(z)$ and $x\notin C(z)$. If $x\in C(z)$, then
$n\le u(z)\le n+N$. Hence
 $z\in U_{u(z)}\subset U_n\cup\cdots\cup U_{n+N}.$
If $x\notin C(z)$, then  letting $u_0=u(z)$ we must have $u_0<n$.
Let $u_1=u(z_{u_0})$.  If $u_0+u_1< n$ we take
$u_2=u(z_{u_0+u_1})$. We proceed in this way until we find the
first $s\le n $ such that $n\le u_0+ \dots +u_s$. Note that
$u_s=u(z_{u_0+\cdots +u_{s-1}})$, and by the choice of $s$ we must
have $x\in C(z_{u_0+\cdots +u_{s-1}})$. Our assumption implies
that $u(z_{u_0+\cdots +u_{s-1}})\le N$, and so $u_0+ \dots +u_s\le
n+N$. By construction we have
%
$$
        z   \in   U_{u_0} $$
        $$
     f^{u_0}(z)=z_{u_0}  \in  U_{u_1} $$
     $$
 f^{u_0+u_1}(z)=z_{u_0+u_1}  \in  U_{u_2}$$
 $$\vdots $$
 $$
 f^{u_0+\cdots u_{s-1}}(z)=z_{u_0+\cdots u_{s-1}}\in U_{u_s}
 $$
By the definition of a  concatenated collection we conclude that
$z\in U_{u_0+u_1+\dots+u_s}$. \cqd

\section{Proofs of main results}\label{se.general}

Let us now prove Theorem \ref{Mtheo}. Suppose that $h\in
L^p(\leb)$, for some $p>3$. This implies that
$\sum_{n\ge1}n^p\leb(h^{-1}(n))<\infty$, and so there exists some
constant $K>0$ such that $$\leb(h^{-1}(n))\le
Kn^{-p},\quad\text{for every $n\ge1$.}$$ Now,  taking
$0<\gamma<(p-3)/(p-1)$ we have for some $K'>0$
 $$
\sum_{n=1}^{\infty}n \left(\sum_{k=
n}^{\infty}\leb(h^{-1}(k))\right)^{1-\gamma}\le
\sum_{n=1}^{\infty}n (K'/n^{p-1})^{1-\gamma} <\infty.
 $$
Defining $$U_n=\{x\in M\ \colon |\det Df^n(x)|\ge b_n\},$$ then we
have that $(U_n )_n$ is a concatenated collection with respect to
the Lebesgue measure. Moreover, setting $$U^*_n= U_n\setminus(
 U_1\cup...\cup U_{n-1})$$ one has
$U^*_n\subset \bigcup_{m\ge n}h^{-1}(m)$, for otherwise there
would be $x\in U^*_n\cap h^{-1}(m)$ with $m<n$, and so $a_m\ge
b_m>|\det Df^m(x)|\ge a_m,$ which is not possible. As $|\det
Df^j(x)|< b_j$ for every $x\in U^*_n$ and $j<n$, we get
$\leb(f^j(U^*_n))\le b_j \leb(U^*_n)$ for each $j<n$. Hence
\begin{align*}
\sum_{n=n_0+1}^{\infty}\sum_{j=0}^{n-1} \leb(f^j(U^*_n))&\le
\sum_{n=n_0+1}^{\infty}\sum_{j=0}^{n-1}b_j \leb(U^*_n)\\
&\le
\sum_{n=n_0+1}^{\infty}\sum_{j=0}^{n_0-1}b_j \leb(U^*_n)+\sum_{n=n_0+1}^{\infty}\sum_{j=n_0}^{n-1}b_j \leb(U^*_n)\\
&\le
\sum_{j=0}^{n_0-1}b_j+\sum_{n=n_0+1}^{\infty}\sum_{j=n_0}^{n-1}b_j
\leb(U^*_n)
\end{align*}
Now we just have to check that the last term in the sum above is
finite. Indeed,
\begin{align*}
\sum_{n=n_0+1}^{\infty}\sum_{j=n_0}^{n-1}b_j \leb(U^*_n) &\le
\sum_{n=n_0+1}^{\infty}\sum_{j=n_0}^{n-1}b_j\sum_{k=
n}^{\infty} \leb(h^{-1}(k))\\
&\le\sum_{n=n_0+1}^{\infty}n b_n\sum_{k= n}^{\infty} \leb(h^{-1}(k))\\
&\le\sum_{n=n_0+1}^{\infty}n \left(\sum_{k=n}^{\infty}
\leb(h^{-1}(k)\right)^{-\gamma}\sum_{k=
n}^{\infty} \leb(h^{-1}(k))\\
&=\sum_{n=n_0+1}^{\infty}n \left(\sum_{k= n}^{\infty}
\leb(h^{-1}(k))\right)^{1-\gamma}<\infty.
\end{align*}
Applying Lemmas~\ref{JClemma1}~and~\ref{JClemma2}, we get for each
generic point $x\in M$ a positive integer number $N_x$ such that
if $y\in f^{-n}(x)$ then $y\in U_{n+s}$ for some $0\le s\le N_x$.
Therefore, $|\det Df^{n+s}(y)|>b_{n+s}\ge b_{n}$. Then, taking
$C_x=K^{-N_x}$, where $K=\sup\{|\det Df(z)|\colon  z\in M\},$ we
obtain the conclusion of Theorem~\ref{JCtheo}:
$$|\det Df^{n}(y)|=\frac{|\det Df^{n+s}(y)|}{|\det
Df^{s}(x)|}>C_x b_{n}.$$ 

Now we explain how we use Theorem~\ref{JCtheo} to prove
Corollary~\ref{Mtheo}. Recall that in Corollary~\ref{Mtheo} we
have $a_n=e^{\lambda n}$ for each $n\in \NN$.
 Assume first that $\leb(\Gamma_n)\le \co(e^{-c^\prime n})$ for
some $c^\prime>0$. Then it is possible to choose $c>0$ such that
 $b_n =e^{cn},$ for $n\ge n_0$.
 The other two  cases
are obtained under similar considerations.


\section{Examples: non-uniformly expanding maps}\label{se.examples}

An important class of dynamical systems where we can immediately
apply our results are the non-uniformly expanding dynamical maps
introduced in \cite{ABV}. As particular examples of this kind of
systems we present below quadratic maps and the higher dimensional
Viana maps.

\paragraph*{Quadratic maps.}

Let $f_a\colon [-1,1]\to [-1,1]$ be given by $f_a(x)=1-ax^2$, for
$0<a\le 2$. Results in \cite{BC1,J} give that for  a positive
Lebesgue measure set of parameters $f_a$ in non-uniformly
expanding. Ongoing work \cite{F} gives that for a positive
Lebesgue measure set of parameters there are $C,c>0$ such that
$\leb(\Gamma_n)\le Ce^{-c n}$ for every $n\ge1$.

Thus, it follows from  Corollary~\ref{Mtheo} that {\em we may find
$\beta>0$ such for Lebesgue almost every $x\in I$ there is $C_x>0$
such that $| (f^n)'(y)|>C_x e^{\beta n}$ for every $y\in
f^{-n}(x)$.}
\paragraph*{Viana maps.}

  Let $a_0\in(1,2)$ be such that the critical point $x=0$
is pre-periodic for the quadratic map $Q(x)=a_0-x^2$. Let
$S^1=\RR/\ZZ$ and $b:S^1\rightarrow \RR$ given by $b(s)=\sin(2\pi
s)$. For fixed small $\alpha>0$, consider the map $\hat f$ from
$S^1\times\RR$ into itself given by $\hat f(s, x) = \big(\hat
g(s),\hat q(s,x)\big)$,
 where  $\hat q(s,x)=a(s)-x^2$ with
$a(s)=a_0+\alpha b(s)$, and $\hat g$ is the uniformly expanding
map of $S^1$ defined by $\hat{g}(s)=ds$ (mod $\ZZ$) for some
integer $d\ge2$. For $\alpha>0$ small enough there is an interval
$I\subset (-2,2)$ for which $\hat f(S^1\times I)$ is contained in
the interior of $S^1\times I$. Thus, any map $f$ sufficiently
close to $\hat f$ in the $C^0$ topology has $S^1\times I$ as a
forward invariant region. Moreover, there are $C,c>0$ such that
$\leb(\Gamma_n)\le Ce^{-c\sqrt n}$ for every $n\ge1$; see
\cite{AA,BST,V}.

Thus, it follows from  Corollary~\ref{Mtheo} that {\em we may find
$\beta>0$ such for Lebesgue almost every $X\in S^1\times I$ there
is a constant $C_X>0$ such that $| \det Df^n(Y)|>C_X e^{\beta\sqrt
n}$ for every $Y\in f^{-n}(X)$.}


\end{document}